\documentclass[12pt,thmsa]{article}

\begin{document}

\begin{center}
{\LARGE Representations of $U_{\epsilon}^{res}(sl_2)$ via Restricted\\[0pt]
q-Fock Spaces}\\[0pt]
\bigskip Xufeng Liu $^{a}$ and Changpu Sun $^{b}$\\[0pt]
\medskip $^{a}$ Department of Mathematics, Peking University\\[0pt]
Beijing 100871, P.R.China\\[0pt]
$^{b}$ Institute of Theoretical Physics, Chinese Academy of Sciences\\[0pt]
Beijing 100080, China
\end{center}

\vspace{20mm}

\begin{center}
\textbf{Abstract}
\end{center}

\begin{quotation}
Two restricted $C[q,q^{-1}]-$forms of the well known q-boson algebra are
introduced and the corresponding restricted q-Fock spaces defined. All of
the irreducible highest weight representations, including the infinite
dimensional ones, of $U_\epsilon ^{res}(sl_2)$ of type 1 are constructed
through the restricted q-Fock spaces.
\end{quotation}

\newpage

\section{Introduction and Notations}

In this paper we take $q$ and $q^{-1}$ to be two indeterminates.We denote by
$C(q)$ the field of rational functions of the indeterminate $q$ and denote
by $C[q,q^{-1}]$ the ring of Laurent polynomials in the indeterminates $q$
and $q^{-1}$.As usual,the integer and nonnegative integer sets are denoted
by $Z$ and $Z^{+}$ respectively, while the positive integer set, or the
natural number set, is denoted by $N.$ For $n\in Z,m\in N$ we use the
following notations

\begin{eqnarray*}
\lbrack n]_q &=&\frac{q^n-q^{-n}}{q-q^{-1}},\ [m]_q!=[m]_q[m-1]_q\cdots
[1]_q, \\
\left[
\begin{array}{c}
n \\
m
\end{array}
\right] _q &=&\frac{[n]_q[n-1]_q\cdots [n-m+1]_q}{[m]_q!}.
\end{eqnarray*}
For $m\in Z\ \backslash \ N$ We adopt the convention
\[
\left[
\begin{array}{c}
n \\
m
\end{array}
\right] _q=0.
\]
For $\varepsilon \in C,$ $\left[
\begin{array}{c}
n \\
m
\end{array}
\right] _\varepsilon $ denotes the complex number obtained from $\left[
\begin{array}{c}
n \\
m
\end{array}
\right] _q$ by substituting $q=\varepsilon $ into the expression, as the
symbol suggests.

Let $U_q$ be an associative algebra over $C(q)$. Naively, it is natural to
think of $U_q$ as a family of algebras depending on a ''parameter''
q.Mathematically, this can be made precise as follows.If $\epsilon \in C$ is
transcendental,one can specialize the indeterminate $q$ to $\epsilon $ by
defining $U_\epsilon =U_q\otimes _{C(q)}C$,via the algebra homomorphism $%
C(q)\longrightarrow C$ that takes $q$ to $\epsilon $. When $\epsilon $ is
algebraic the above homomorphism from $C(q)$ to $C$ is not available and
this direct specialization of $U_q$ might not make sense. Nevertheless, one
can proceed by first constructing a $C[q,q^{-1}]$--form ,or integral form,
of $U_q$ ,namely, a $C[q,q^{-1}]$--subalgebra $\tilde{U}_q$ of $U_q$ such
that $U_q=\tilde{U}_q\otimes _{C[q,q^{-1}]}C(q)$.Then one defines the
specialization $U_\epsilon $ of $U_q$ as $\tilde{U}_q\otimes _{C[q,q^{-1}]}C$%
,via the algebra homomorphism $C[q,q^{-1}]\longrightarrow C$ that takes $q$
to $\epsilon $.

The quantum algebra $U_q(g)$ associated to a Kac-Moody algebra $g$ is an
associative algebra over $C(q).$When $\epsilon $ is not a root of unity the
representation theory of $U_\epsilon (g)$ has been well established[1]. To
deal with the case where $\epsilon $ is a root of unity two $C[q,q^{-1}]$%
--forms of $U_q(g)$,namely,the ``non-restricted'' form and the
``restricted'' form,have been introduced.In the ``non-restricted'' form,one
takes $\tilde{U}_q(g)$ to be the $C[q,q^{-1}]$--subalgebra of $U_q(g)$
generated by the Chevalley generators $e_i,f_i$ and some other elements of $%
U_q(g)$.In this case,the finite dimensional representations of $U_\epsilon
(g)$ have been studied by De Concini,Kac and Procesi,when $g$ is finite
dimensional,and by Beck and Kac,when $g$ is untwisted affine[2,3]. The
restricted $C[q,q^{-1}]-$form $U_q^{res}(g)$ of $U_q(g)$ is introduced by
Lusztig [4]. The study of representation theory of $U_\epsilon ^{res}(g)$ is
pioneered also by Lusztig and developed by Chari and Pressley [5,6,7].

The q-boson realization method has been widely used to construct
representations of quantum algebras in both generic case and root of unity
case [8-16]. Especially, cyclic q-boson algebra has been introduced to
obtain the so called cyclic representations of $U_\epsilon (g)$ associated
with the non-restricted $C[q,q^{-1}]-$form of $U_q(g)\ [15,16].$ In this
paper, we will introduce two restricted $C[q,q^{-1}]-$forms of the well
known q-boson algebra and define the restricted q-Fock spaces
correspondingly. Then we will construct all irreducible highest weight
representations of $U_\epsilon (sl_2)$ of type 1 on the restricted q-Fock
spaces.

\section{Some Basic Facts about $U_{\mathcal{\epsilon }}^{res}(sl_2)$}

Let us recall some definitions and basic facts. For details we refer readers
to Ref.[1].Throughout this paper we use the notation $\mathcal{A}=C[q,q^{-1}]
$.

\textbf{Definition 2.1.} The quantum algebra $U_q(sl_2)$ is the associative
algebra over $C(q)$ with generators $e,f,K$ and $K^{-1}$ and the following
relations:
\begin{eqnarray*}
&&KK^{-1}=K^{-1}K=1, \\
&&KeK^{-1}=q^2e,KfK^{-1}=q^{-2}f, \\
&&\left[ e,f\right] =\frac{K-K^{-1}}{q-q^{-1}}.
\end{eqnarray*}

\textbf{Definition 2.2.} The algebra $U_{\mathcal{A}}^{res}(sl_2)$ is the $%
\mathcal{A}$-subalgebra of $U_q(sl_2)$ generated by the elements $%
e^{(r)},f^{(r)},K^{\pm 1}$ $(r\in N)$ for $r\geq 1$, where $e^{(r)}=\frac{e^r%
}{[r]_q!}$ and $f^{(r)}=\frac{f^r}{[r]_q!}.$

From now on we assume that $\epsilon $ is a primitive$\mathcal{\ }p$th root
of unity, where $p$ is odd and greater than 1. When necessary, $C$ is
considered as $\mathcal{A-}$module via the the algebra homomorphism $%
\mathcal{A}\longrightarrow C$ that takes $q$ to $\epsilon .$By definition,
the restricted specialization of $U_q(sl_2)$ is
\[
U_{\mathcal{\epsilon }}^{res}(sl_2)=U_{\mathcal{A}}^{res}(sl_2)\otimes _{%
\mathcal{A}}C.
\]
For simplicity, $e^{(r)}\otimes _{\mathcal{A}}1,f^{(r)}\otimes _{\mathcal{A}%
}1$ and $K^{\pm 1}\otimes _{\mathcal{A}}1$ will be identified with $%
e^{(r)},f^{(r)}$ and $K^{\pm 1}$ respectively.

If $V$ is a representation of $U_{\mathcal{\epsilon }}^{res}(sl_2)$ on which
$K$ is diagonalizable and $K^{\pm p}=1$ it is said to be of type $1.$

\textbf{Definition 2.3.} Let $m$ be an integer. The weight space $V_m$ of a $%
U_{\mathcal{\epsilon }}^{res}(sl_2)-$module $V$ of type $1$ is defined by
\[
V_m=\{v\in V|\ Kv=\epsilon ^mv,\ \left[
\begin{array}{c}
K;0 \\
p
\end{array}
\right] _qv=\left[
\begin{array}{c}
m \\
p
\end{array}
\right] _\epsilon v\},
\]
where
\[
\left[
\begin{array}{c}
K;0 \\
p
\end{array}
\right] _q=\prod_{s=1}^p\frac{Kq^{1-s}-K^{-1}q^{s-1}}{q^s-q^{-s}}
\]
belongs to $U_{\mathcal{A}}^{res}(sl_2)$ and is identified with $\left[
\begin{array}{c}
K;0 \\
p
\end{array}
\right] _q\otimes _{\mathcal{A}}1.$

For any integer $n$, write $n=n_0+pn_1,$ where $n_0$ and $n_1$ are integers
and $0\leq n_0<p.$ It is readily verified that
\[
\left[
\begin{array}{c}
n \\
p
\end{array}
\right] _\epsilon =n_1.
\]
Then by definition, $v\in V$ is a weight vector of weight $m$ and we have
\[
Kv=\epsilon ^{m_0}v,\ \left[
\begin{array}{c}
K;0 \\
p
\end{array}
\right] _qv=m_1v.
\]

A $U_{\mathcal{\epsilon }}^{res}(sl_2)-$module $V$ of type $1$ is called a
highest weight module if it is generated by a primitive vector, i.e. a
vector $v_\lambda $ for some $\lambda \in Z$ such that $ev_\lambda
=e^{(p)}v_\lambda =0.$ It is obvious that for such a module $V$ we have
\[
V=\sum_{\mu \leq \lambda }\oplus \ V_\mu
\]
so that $\lambda $ is the highest weight of $V$ and $V$ $_\lambda
=Cv_\lambda .$ It then follows by the usual argument that $V$ has a unique
irreducible quotient module.

Denote by $V_\epsilon ^{res}(\lambda )$ the irreducible highest weight $U_{%
\mathcal{\epsilon }}^{res}(sl_2)-$module of type $1$ and of highest weight $%
\lambda .$ Then by the Verma module construction and the above argument one
can prove that $V_\epsilon ^{res}(\lambda )$ is isomorphic to $V_\epsilon
^{res}(\mu )$ if $\lambda =\mu .$

\section{Restricted q-Fock Spaces}

The $q-$boson algebra $B_q(n)$of rank $n$ is the associative algebra over $%
C(q)$ generated by the elements $a_i,a_i^{+},K_i^{\pm 1},1$ $(i=1,2,\cdots
,n)$ with the following relations
\begin{eqnarray*}
a_ia_i^{+}-q^2a_i^{+}a_i &=&1,[a_i,a_j^{+}]=0\,(i\neq j) \\
\lbrack a_i,a_j] &=&[a_i^{+},a_j^{+}]=0 \\
K_ia_j^{+}K_i^{-1} &=&q^{\delta _{ij}}a_j^{+},K_ia_jK_i^{-1}=q^{-\delta
_{ij}}a_j \\
K_iK_i^{-1} &=&K_i^{-1}K_i=1
\end{eqnarray*}
where $1$ is the unit.

We introduce the notations
\[
a_i^{(n)}=\frac{a_i^n}{[n]_q!},\,a_i^{+(n)}=\frac{a_i^{+n}}{[n]_q!}
\]
and the following definition.

\textbf{Definition 3.1.} The type 1 restricted $q-Boson$ algebra $%
B_q^{res1}(2)$ of rank 2 is the $\mathcal{A-}$subalgebra of $B_q(2)$
generated by the elements $a_i,a_i^{+(r)},K_i^{\pm 1},1$ $(i=1,2;r\in N);$%
The type 2 restricted $q-Boson$ algebra $B_q^{res2}(2)$ of rank 2 is the $%
\mathcal{A-}$subalgebra of $B_q(2)$ generated by the elements $%
a_1,a_1^{+(r)},a_2^{(r)},a_2^{+},K_1^{\pm 1},K_2^{\pm 1},1$ $(r\in N).$

By induction one can prove the following three lemmas without difficulty.

\textbf{Lemma 3.1.}
\[
\left[
\begin{array}{l}
r \\
k
\end{array}
\right] _q=q^{-k}\left[
\begin{array}{c}
r-1 \\
k
\end{array}
\right] _q+q^{r-k}\left[
\begin{array}{c}
r-1 \\
k-1
\end{array}
\right] _qfor\ r\geq k\geq 0.
\]

\textbf{Lemma 3.2.}
\[
\left[
\begin{array}{l}
r \\
k
\end{array}
\right] _q\in C[q,q^{-1}]\,for\,r\geq k\geq 0.
\]

\textbf{Lemma 3.3.}
\begin{eqnarray*}
a_i^na_i^{+(m)} &=&\sum_{s=0}^mq^{n(s+m)}q^{\frac{(s-m)(s+m+1)}2}\left[
\begin{array}{c}
n \\
m-s
\end{array}
\right] _qa^{+(s)}a^{n-m+s} \\
a_i^{(n)}a_i^{+m} &=&\sum_{s=0}^mq^{n(s+m)}q^{\frac{(s-m)(s+m+1)}2}\left[
\begin{array}{c}
m \\
s
\end{array}
\right] _qa^{+s}a^{(n-m+s)}
\end{eqnarray*}

\textbf{Proposition 3.1.} $%
\{a_1^{+(r_1)}a_2^{+(r_2)}a_1^{s_1}a_2^{s_2}K_1^{t_1}K_2^{t_2}|r_i,s_i\in
Z^{+},t_i\in Z;i=1,2\}$ is an $\mathcal{A}-$basis of $B_q^{res1}(2);%
\{a_1^{+(r_1)}a_2^{+r_2}a_1^{s_1}a_2^{(s_2)}K_1^{t_1}K_2^{t_2}|r_i,s_i\in
Z^{+},t_i\in Z;i=1,2\}$ is an $\mathcal{A}-$basis of $B_q^{res2}(2).$

Proof. It follows directly from Lemma 3.2 and Lemma 3.3.

\textbf{Corollary}. $B_q^{res1}(2)$ and $B_q^{res2}(2)$ are two integral
forms of $B_q(2).$

Consider $B_q^{res1}(2)$ as $B_q^{res1}(2)-$module as well as $C[q,q^{-1}]-$%
algebra .Let $I$ be the left ideal of $B_q^{res1}(2)$ generated by the
elements
\[
K_1-1,K_2-1,a_1,a_2.
\]
Obviously, $I$ is $B_q^{res1}(2)-$submodule. Let $\mathcal{F}_{\mathcal{A}%
}^1(2)$ denote the quotient module $B_q^{res1}(2)/I.$ It follows from
Proposition 3.1 that $\mathcal{F}_{\mathcal{A}}^1(2)$ is a free $\mathcal{A}-
$module and

\[
\mathcal{F}_{\mathcal{A}}^1(2)=span_{\mathcal{A}}\{a_1^{+(r_1)}a_2^{+(r_2)}%
\left| 0\right\rangle \,|\,a_1\left| 0\right\rangle =a_2\left|
0\right\rangle =0,K_1\left| 0\right\rangle =K_2\left| 0\right\rangle =\left|
0\right\rangle ;r_1,r_2\in Z^{+}\}.
\]
where $\{a_1^{+(r_1)}a_2^{+(r_2)}\left| 0\right\rangle |r_1,r_2\in Z^{+}\}$
is a basis of $\mathcal{F}_{\mathcal{A}}^1(2).\ \mathcal{F}_{\mathcal{A}%
}^1(2)$ is called $\mathcal{A}$-Fock module associated with $B_q^{res1}(2).\
$The $\mathcal{A}$-Fock module $\mathcal{F}_{\mathcal{A}}^2(2)$ associated
with $B_q^{res2}(2)$ can be defined in a similar way :
\[
\mathcal{F}_{\mathcal{A}}^2(2)=span_{\mathcal{A}}\{a_1^{+(r_1)}a_2^{+r_2}%
\left| 0\right\rangle \,|\,a_1^r\left| 0\right\rangle =a_2^{(r)}\left|
0\right\rangle =0,r\in N;K_1\left| 0\right\rangle =K_2\left| 0\right\rangle
=\left| 0\right\rangle ;r_1,r_2\in Z^{+}\}.
\]

It is easy to see that $\mathcal{F}_{\mathcal{A}}^1(2)\otimes _{\mathcal{A}%
}C(q)=\mathcal{F}_{\mathcal{A}}^2(2)\otimes _{\mathcal{A}}C(q)\stackrel{%
\triangle}{=} \mathcal{F}_q(2)$ is the ordinary q-Fock space associated with
$B_q(2).\ $It is a vector space over the field $C(q).\ $If we regard $%
\mathcal{F}_q(2)$ as $\mathcal{A}-$module,then by identifying $\mathcal{F}_{%
\mathcal{A}}^1(2)$ with $\mathcal{F}_{\mathcal{A}}^1(2)\otimes _{\mathcal{A}%
}1$ , $\mathcal{F}_{\mathcal{A}}^1(2)$ becomes a submodule of $\mathcal{F}%
_q(2),$ and $\mathcal{F}_{\mathcal{A}}^2(2)$ becomes a submodule of $%
\mathcal{F}_{\mathcal{A}}^2(2)$ likewise. For convenience, in the subsequent
discussion we will use the following notations:
\begin{eqnarray*}
f(r_1,r_2) &=&a_1^{+(r_1)}a_2^{+(r_2)}\left| 0\right\rangle , \\
g(r_1,r_2) &=&a_1^{+(r_1)}a_2^{+r_2}\left| 0\right\rangle .
\end{eqnarray*}

We define
\begin{eqnarray*}
\mathcal{F}_{\mathcal{\epsilon }}^1(2) &=&\mathcal{F}_{\mathcal{A}%
}^1(2)\otimes _{\mathcal{A}}C, \\
\mathcal{F}_{\mathcal{\epsilon }}^2(2) &=&\mathcal{F}_{\mathcal{A}%
}^2(2)\otimes _{\mathcal{A}}C,
\end{eqnarray*}
and identify $f(r_1,r_2)\otimes _{\mathcal{A}}1$ and $g(r_1,r_2)\otimes _{%
\mathcal{A}}1$with $f(r_1,r_2)$ and $g(r_1,r_2)$ respectively. $\mathcal{F}_{%
\mathcal{\epsilon }}^1(2)$ and $\mathcal{F}_{\mathcal{\epsilon }}^2(2)$ will
be referred to as restricted q-Fock spaces. In the following sections, we
will construct irreducible representations of $U_{\mathcal{\epsilon }%
}^{res}(sl_2)$ on the restricted q-Fock spaces.

\section{Finite Dimensional Representations of $U_{\mathcal{\epsilon }%
}^{res}(sl_2)$ via Restricted q-Fock Space}

We recall that $U_q(sl_2)$ has the following realization on $\mathcal{F}%
_q(2)\ [8]:$%
\begin{eqnarray*}
e &=&K_2^{-1}a_1^{+}a_2,\,f=K_1^{-1}a_1a_2^{+}, \\
K &=&K_1K_2^{-1},\,K^{-1}=K_1^{-1}K_2.
\end{eqnarray*}

\textbf{Remark. }The above realization looks different from that presented
in Ref.[8]. This is because the generators of the q-boson algebra adopted in
this paper are slightly different from those in Ref.[8].

By the natural action $\mathcal{F}_q(2)$ becomes a $U_q(sl_2)-$module
through this realization.Then naturally it becomes a $U_{\mathcal{A}%
}^{res}(sl_2)-$module through the realization
\[
e^{(r)}=q^{-\frac{r(r-1)}2}K_2^{-r}a_1^{+(r)}a_2^r,\,f^{(r)}=q^{-\frac{r(r-1)%
}2}K_1^{-r}a_2^{+(r)}a_1^r.
\]
It turns out that $\mathcal{F}_{\mathcal{A}}^1(2)$ is a $U_{\mathcal{A}%
}^{res}(sl_2)-$submodule. In fact, the action of $U_{\mathcal{A}}^{res}(sl_2)
$ on $\mathcal{F}_{\mathcal{A}}^1(2)$ reads as follows:
\begin{eqnarray*}
Kf(r_1,r_2)
&=&q^{r_1-r_2}f(r_1,r_2),\,K^{-1}f(r_1,r_2)=q^{r_2-r_1}f(r_1,r_2), \\
e^{(r)}f(r_1,r_2) &=&\left[
\begin{array}{c}
r+r_1 \\
r
\end{array}
\right] _qf(r_1+r,r_2-r), \\
\,f^{(r)}f(r_1,r_2) &=&\left[
\begin{array}{c}
r+r_2 \\
r
\end{array}
\right] _qf(r_1-r,r_2+r).
\end{eqnarray*}
Here we have used Lemma 3.3 and the formula
\[
a_i^{+(m)}a_i^{+(n)}=\left[
\begin{array}{c}
m+n \\
m
\end{array}
\right] _qa_i^{+(m+n)}\,for\,i=1,2;m,n\in N.
\]
Thus it follows that $U_{\mathcal{A}}^{res}(sl_2)\mathcal{F}_{\mathcal{A}%
}^1(2)\subset \mathcal{F}_{\mathcal{A}}^1(2).$ This proves the claim.

Clearly, $\mathcal{F}_{\mathcal{\epsilon }}^1(2)$ is a $U_{\mathcal{A}%
}^{res}(sl_2)-$module. It then follows that it is also a $U_{\mathcal{%
\epsilon }}^{res}(sl_2)-$module. Here the action of $U_{\mathcal{\epsilon }%
}^{res}(sl_2)$ on $\mathcal{F}_{\mathcal{\epsilon }}^1(2)$ is induced from
that of $U_{\mathcal{A}}^{res}(sl_2)$ by identifying $q$ with $\epsilon .$
Explicitly, we have
\begin{eqnarray*}
Kf(r_1,r_2) &=&\epsilon ^{r_1-r_2}f(r_1,r_2),\,K^{-1}f(r_1,r_2)=\epsilon
^{r_2-r_1}f(r_1,r_2), \\
e^{(r)}f(r_1,r_2) &=&\left[
\begin{array}{c}
r+r_1 \\
r
\end{array}
\right] _\epsilon f(r_1+r,r_2-r), \\
\,f^{(r)}f(r_1,r_2) &=&\left[
\begin{array}{c}
r+r_2 \\
r
\end{array}
\right] _\epsilon f(r_1-r,r_2+r).
\end{eqnarray*}

We observe that for $m\in N$ the subspace
\[
V_m\stackrel{\triangle}{=} span_C\ \{f(r_1,r_2)\ |\ r_1+r_2=m\}
\]
of $\mathcal{F}_{\mathcal{\epsilon }}^1(2)$ is a $U_{\mathcal{\epsilon }%
}^{res}(sl_2)-$submodule.Define $v_i^{(m)}=f(m-i,i)$ for $i=0,1,\cdots ,m.$
Then $\{v_0^{(m)},\cdots ,v_m^{(m)}\}$ is a basis of $V_m.$

The action of $U_{\mathcal{\epsilon }}^{res}(sl_2)$ on $V_m$ is given by
\begin{eqnarray*}
Kv_r^{(m)} &=&\epsilon ^{m-2r}v_r^{(m)},\ K^{-1}v_r^{(m)}=\epsilon
^{2r-m}v_r^{(m)}, \\
ev_r^{(m)} &=&\left[ m-r+1\right] _\epsilon v_{r-1}^{(m)},\
fv_r^{(m)}=\left[ r+1\right] _\epsilon v_{r+1}^{(m)}, \\
e^{(p)}v_r^{(m)} &=&((m-r)_1+1)v_{r-p}^{(m)},\ f^{(p)}v_r^{(m)}=\left(
r_1+1\right) \ v_{r+p}^{(m)}.
\end{eqnarray*}
Here $v_i^{(m)}$ is understood as the zero vector if $i<0$ or $i>m.$ We
notice that $V_m$ is none other than the so called Weyl module $W_\epsilon
^{res}(m)$ with maximal weight $m.$ Consequently, we have the following
result.

\textbf{Proposition 4.1}. $V_m$ is irreducible if and only if either $m<p$
or $m_0=p-1;$ If $V_m$ is reducible it is not completely reducible.

A sketch of the proof of this proposition can be found in Ref.[1].The main
point is the observation that
\[
V^{\prime }\stackrel{\triangle }{=}span_C\{v_r^{(m)}\ |\ m_0<r_0<p,\
r_1<m_1\}
\]
is the unique proper $U_{\mathcal{\epsilon }}^{res}(sl_2)-$submodule of $V_m.
$

It is an established fact that every finite dimensional irreducible $U_{%
\mathcal{\epsilon }}^{res}(sl_2)-$module of type $1$ is a quotient module of
Weyl module. As we have realized the Weyl module $W_\epsilon ^{res}(m)$ for
an arbitrary $m\in N$ via the q-Fock space, we have actually presented
q-Fock space construction for all finite dimensional irreducible
representations of type $1$of $U_{\mathcal{\epsilon }}^{res}(sl_2).$

\section{Infinite Dimensional Representations of $U_{\mathcal{\epsilon }%
}^{res}(sl_2)$ via Restricted q-Fock Space}

In this section let us turn to consider q-Fock space construction for some
infinite dimensional irreducible representations of $U_{\mathcal{\epsilon }%
}^{res}(sl_2).$ To this end, we need to use another q-boson realization of $%
U_q(sl_2).$ Our starting point is the following well-known q-boson
realization of $U_q(sl_2)$ on $\mathcal{F}_q(2)$ [13] :
\begin{eqnarray*}
K &=&q^{-1}K_1^{-1}K_2^{-1},\ K^{-1}=qK_1K_2, \\
e &=&K_1^{-1}K_2^{-1}a_1a_2,\ f=-a_1^{+}a_2^{+}.
\end{eqnarray*}
By definition $\mathcal{F}_q(2)$ is a $U_q(sl_2)-$module, and thus a $U_{%
\mathcal{A}}^{res}(sl_2)-$module, the action on $\mathcal{F}_q(2)$ being the
natural one.We have
\[
e^{(r)}=q^{-r(r-1)}K_1^{-r}K_2^{-r}a_1^ra_2^{(r)},\
f^{(r)}=(-1)^ra_1^{+(r)}a_2^{+r}.
\]
Then it follows from Lemma 3.3 that
\begin{eqnarray*}
e^{(r)}g(r_1,r_2) &=&\left[
\begin{array}{c}
r_2 \\
r_2-r
\end{array}
\right] _qg(r_1-r,r_2-r), \\
f^{(r)}g(r_1,r_2) &=&(-1)^r\left[
\begin{array}{c}
r+r_1 \\
r
\end{array}
\right] _qg(r_1+r,r_2+r),
\end{eqnarray*}
namely, $U_{\mathcal{A}}^{res}(sl_2)\mathcal{F}_{\mathcal{A}}^2(2)\subset
\mathcal{F}_{\mathcal{A}}^2(2).$ Consequently, $\mathcal{F}_{\mathcal{A}%
}^2(2),$ and hence $\mathcal{F}_{\mathcal{\epsilon }}^2(2)$ is a $U_{%
\mathcal{A}}^{res}(sl_2)-$module. Finally, $\mathcal{F}_{\mathcal{\epsilon }%
}^2(2)$ becomes a $U_{\mathcal{\epsilon }}^{res}(sl_2)-$module in the
obvious way. Explicitly, the action of $U_{\mathcal{\epsilon }}^{res}(sl_2)$
on $\mathcal{F}_{\mathcal{\epsilon }}^2(2)$ is as follows:
\begin{eqnarray*}
Kg(r_1,r_2) &=&\epsilon ^{-(r_1+r_2+1)}g(r_1,r_2),\
K^{-1}g(r_1,r_2)=\epsilon ^{(r_1+r_2+1)}g(r_1,r_2), \\
e^{(r)}g(r_1,r_2) &=&\left[
\begin{array}{c}
r_2 \\
r_2-r
\end{array}
\right] _\epsilon g(r_1-r,r_2-r), \\
f^{(r)}g(r_1,r_2) &=&(-1)^r\left[
\begin{array}{c}
r+r_1 \\
r
\end{array}
\right] _\epsilon g(r_1+r,r_2+r).
\end{eqnarray*}

On $\mathcal{F}_{\mathcal{\epsilon }}^2(2),$ we have
\[
\left[
\begin{array}{c}
K;0 \\
p
\end{array}
\right] _q=\prod_{s=1}^p\frac{Kq^{1-s}-K^{-1}q^{s-1}}{q^s-q^{-s}}%
=\prod_{s=1}^p\frac{K_1^{-1}K_2^{-1}q^{-s}-K_1K_2q^s}{q^s-q^{-s}}.
\]
Thus
\[
\left[
\begin{array}{c}
K;0 \\
p
\end{array}
\right] _qg(r_1,r_2)=\left[
\begin{array}{c}
-r_1-r_2-1 \\
p
\end{array}
\right] _\epsilon g(r_1,r_2)=-(r_1+r_2+1)_1g(r_1,r_2)
\]
Hence, $g(r_1,r_2)$ is a weight vector of the weight $-(r_1+r_2+1)$ and $%
\mathcal{F}_{\mathcal{\epsilon }}^2(2)$ is a $U_{\mathcal{\epsilon }%
}^{res}(sl_2)-$module of type $1.$

For an arbitrary integer $s,$ define
\[
V^s=span_C\{g(r_1,r_1+s)|r_1\in Z^{+}\}.
\]
We observe that $V^s$ is a $U_{\mathcal{\epsilon }}^{res}(sl_2)-$submodule
and there is the following decomposition:
\[
\mathcal{F}_{\mathcal{\epsilon }}^2(2)=\sum_{s\in Z}\oplus V^s
\]
For $s$ that satisfies $s_0=0,\ $namely, $s=\pm rp\ (r\in Z^{+}),$ we have
the following result.

\textbf{Proposition 5.1.} For $r\in Z^{+},$ $V^{\pm rp}$ is isomorphic to $%
V_\epsilon ^{res}(-(rp+1)).$

Proof. We only need to consider the case of $V^{rp}.$ The other case is
similar. We first show that $V^{rp}$ is an infinite dimensional irreducible $%
U_{\mathcal{\epsilon }}^{res}(sl_2)-$module.For $m\in Z^{+}$ define $%
g(m)=g(m,m+rp).$ Then $\{g(m)|\ m\in Z^{+}\}$ is a basis of $V^{rp}.$ We
have
\[
eg(m)=[m]_\epsilon g(m-1),\ fg(m)=-[m+1]_qg(m+1)
\]
and
\[
e^{(p)}g(np)=ng((n-1)p),\ f^{(p)}g(np-1)=(n-1)g((n+1)p-1)
\]
for $n\in N.$ Using these relations, one can easily show that as $U_{%
\mathcal{\epsilon }}^{res}(sl_2)-$module $V^{rp}$ can be generated by an
arbitrary vector in it. The irreducibility is thus proved.

As shown above, $g(0)$ is a weight vector of the weight $-(1+rp).$ It is
easy to see that it is actually a highest weight vector. On the other hand,
it is clear that $V^{rp}$ is a highest weight module: $V^{rp}=$ $U_{\mathcal{%
\epsilon }}^{res}(sl_2)g(0).$ The claim then follows from the irreducibility
of $V^{rp}.$

Now let us study $V^s$ with $s_0\neq 0.$ When $s>0,$ define $v_m=g(m,m+s)$
and the subspace $V^{\prime }$ of $V^s:$%
\[
V^{\prime }=span_C\{v_m|\ p>m_0\geq p-s_0\}.
\]

When $s<0,$ we use $W^s$ to denote $V^s$ for convenience. Then
\[
W^s=span_C\{g(r_1+|s|,r_1)|r_1\in Z^{+}\}.
\]
We define $w_m=g(m+|s|,m)$ and the subspace $W^{\prime }$ of $W^s:$%
\[
W^{\prime }=span_C\{w_m|\ 0\leq m_0<p-s_0\}.
\]
It is not difficult to see that $V^{\prime }$ and $W^{\prime }$ are proper
submodules of $V^s$ and $W^s$ respectively.

\textbf{Proposition 5.2.} $V^{\prime }$ is isomorphic to $V_\epsilon
^{res}(\lambda )$ with $\lambda =-(p-s_0+(s_1+1)p+1).$

Proof. By the explicit action of $U_{\mathcal{\epsilon }}^{res}(sl_2)$ on $%
\mathcal{F}_{\mathcal{\epsilon }}^2(2),$ it is not difficult to see that $%
V^{\prime }$ is irreducible and $V^{\prime }=U_{\mathcal{\epsilon }%
}^{res}(sl_2)v_{p-s_0}.$ Since $v_{p-s_0}$ is the highest weight vector of $%
V^{\prime }$ of the weight $-(p-s_0+(s_1+1)p+1),$ the proposition then
follows.

\textbf{Proposition 5.3.} $V^{\prime }$ is the unique maximal proper $U_{%
\mathcal{\epsilon }}^{res}(sl_2)-$submodule of $V^s$ with $s>0$ and $s_0\neq
0.$

Proof. Let $V^{\prime \prime }$ be an arbitrary proper $U_{\mathcal{\epsilon
}}^{res}(sl_2)-$submodule of $V^s.$ To prove the proposition, we only need
to show that it is included in $V^{\prime }.$We observe that $V^{\prime
\prime }$ has the vector space decomposition
\[
V^{\prime \prime }=\sum_{\lambda \in \Lambda _s}\oplus \ (V_\lambda
^s\bigcap V^{\prime \prime })
\]
where $\Lambda _s$ stands for the weight set of $V^s.$ Suppose $v$ is a
vector in $V^{\prime \prime }.$ Define
\begin{eqnarray*}
\Lambda _s^{\prime } &=&\{\lambda \in \Lambda _s|\ \lambda _0\geq p-s_0\}, \\
\Lambda _s^{\prime \prime } &=&\{\lambda \in \Lambda _s|\ \lambda _0<p-s_0\}.
\end{eqnarray*}
Then we can write $v=v_1+v_2$ where
\[
v_1\in \sum_{\lambda \in \Lambda _s^{\prime \prime }}\oplus \ V_\lambda
^s,\,v_2\in \sum_{\lambda \in \Lambda _s^{\prime }}\oplus \ V_\lambda ^s\in
V^{\prime }.
\]
Now due to the above vector space decomposition of $V^{\prime \prime },$ we
have $v_1\in V^{\prime \prime }.$ On the other hand, if $v_1$ is a nonzero
vector, it is not difficult to show that $U_{\mathcal{\epsilon }%
}^{res}(sl_2)v_1=V^s.$ Consequently, $v_1=0$ as $V^{\prime \prime }$ is a
proper $U_{\mathcal{\epsilon }}^{res}(sl_2)-$submodule by assumption. This
means $v=v_2\in V^{\prime },$ namely, $V^{\prime \prime }\subset V^{\prime
}. $ The proposition is thus proved.

\textbf{Corollary.} For $V^s$ with $s>0$ and $s_0\neq 0,$ $V^s/V^{\prime }$
is isomorphic to $V_\epsilon ^{res}(-(s+1)).$

Proof. Since $V^{\prime }$ is a maximal proper submodule, $V^s/V^{\prime }$
is irreducible. Moreover, $V^s/V^{\prime }$ is a highest module with the
highest vector $\overline{v}_0\stackrel{\triangle}{=} v_0+V^{\prime
}=g(0,s)+V^{\prime }$ of the weight $-(s+1).$ The corollary thus follows.

\textbf{Remark.} Obviously, the above studied $U_{\mathcal{\epsilon }%
}^{res}(sl_2)-$modules $V^{np},V^{\prime }$ and $V^s/V^{\prime }$ are all
infinite dimensional. We notice that they are isomorphic to the modules $%
V_\epsilon ^{res}(\lambda )$ with $\lambda $ negative. This should be the
case. Actually, if $\lambda $ is a positive integer $V_\epsilon
^{res}(\lambda )$ must be finite dimensional.

In the same way, we can prove the following result.

\textbf{Proposition 5.4.} $W^{\prime }$ is isomorphic to $V_\epsilon
^{res}(-(|s|+1));\ W^s/W^{\prime }$ is isomorphic to $V_\epsilon
^{res}(\lambda )$ with $\lambda =-(p-|s|_0+(1+|s|_1)p+1).$

In summary, we have constructed all the representations $V_\epsilon
^{res}(\lambda )$ of $U_{\mathcal{\epsilon }}^{res}(sl_2)$ with $\lambda \in
Z$ via the restricted q-Fock spaces. According to the representation theory
of $U_{\mathcal{\epsilon }}^{res}(sl_2)$ they have exhausted all the
irreducible highest representations of $U_{\mathcal{\epsilon }}^{res}(sl_2)$
of type 1.Obviously, the method presented in this paper can readily be
generalized to be applicable to the case of $U_{\mathcal{\epsilon }%
}^{res}(sl_n).$ \newpage References

\begin{enumerate}
\item  V.Chari and A.N.Pressley,A Guide to Quantum Groups,Cambridge
University Press,Cambridge,1994.

\item  C.De Concini and V.G.Kac,Representations of quantum groups at roots
of unity,in Operator Algebras,Unitary Representations,Enveloping Algebras
and Invariant Theory,A.Connes,M.Euflo,A.Joseph and R.Rentschler
(eds.),pp.471-506,Progress in Mathematics 92,Birkhauser,Boston,1990.

\item  J.Beck and V.G.Kac,Finite dimensional representations of quantum
affine algebras at roots of unity,J.Amer.Math.Soc.9(1996)391-423.

\item  G.Lusztig, Quantum deformations of certain simple modules over
enveloping algebras, Adv. Math. 70(1988)237-249.

\item  G.Lusztig, Finite dimensional Hopf algebras arising from quantized
universal enveloping algebras, J. Amer. Math. Soc. 3(1990)257-296.

\item  G.Lusztig, Quantum groups at roots of 1, Geom. Dedicata
35(1990)89-113.

\item  V.Chari and A.N.Pressley, Quantum affine algebras at roots of unity,
Representation Theory 1(1997)280--328.

\item  C.P.Sun and H.C.Fu, The q-deformed boson realization of the quantum
group $SU_q(n)$ and its representations, J.Phys.A 22(1989)L983-986.

\item  L.C.Biedenharn, The quantum group $SU_q(2)$ and a q-analogue of the
boson operators, J.Phys. A 22(1989)L873-878.

\item  A.J.Macfarlane, On q-analogue of the quantum harmonic oscillator and
the quantum group $SU_q(2)$. J.Phys.A 22(1989)4581-4588.

\item  C.P.Sun and M.L.Ge, The q-analogue of the boson algebra, its
representation on the Fock space, and applications to the quantum group,J.
Math. Phys. 32(1991)595-598.

\item  C.P.Sun and M.L.Ge, The q-deformed boson realization of
representations of quantized universal enveloping algebras for q a root of
unity: I.The case of $U_q(SL(l))$, J.Phys. A 24(1991)3265-3280.

\item  C.P.Sun, X.F.Liu and M.L.Ge, A new q-deformed boson realization of
quantum algebra $sl_q(2)$ and nongeneric $sl_q(2)$ R-matrices, J. Math.
Phys. 32(1991)2409-2412.

\item  C.P.Sun and M.L.Ge, The q-deformed boson realization of
representations of quantized universal enveloping algebras for q a root of
unity: II. The subalgebra chain of $U_q(C_l)$, J.Phys. A 25(1992)401-410.

\item  C.P.Sun, X.F.Liu and M.L.Ge, Cyclic boson algebra and q-boson
realization of cyclic representation of quantum algebra $sl_q(3)$ at $q^p=1$%
J.Phys.A 25(1992)161-168.

\item  C.P.Sun and M.L.Ge, Cyclic boson operators and new representations of
$sl_q(2)$ at q a root of unity, J.Phys.A 24(1991)L969-973.
\end{enumerate}

\end{document}